\documentclass{svmult}

\usepackage{makeidx}         
\usepackage{graphicx}        
\usepackage{multicol}        
\usepackage[bottom]{footmisc}
\usepackage{natbib}
\usepackage{url}
\usepackage{multirow}

\renewcommand{\PackageWarningNoLine}[2]{}
\usepackage{amsmath}
\usepackage{amsfonts}
\usepackage{amssymb}
\usepackage{dsfont}
\usepackage{bm}
\usepackage{subfigure}
\usepackage[algo2e,linesnumbered,ruled]{algorithm2e}
\newtheorem{algorithm}[theorem]{Algorithm}

\begin{document}

\title*{Multigrid Preconditioner for Nonconforming Discretization of Elliptic Problems with Jump Coefficients}
\titlerunning{Multigrid for Nonconforming Discretization}
\author{Blanca Ayuso De Dios\inst{1}\and
Michael Holst\inst{2}\and Yunrong Zhu\inst{2} \and
Ludmil Zikatanov\inst{3}
}
\authorrunning{B. Ayuso de Dios \and M. Holst \and Y. Zhu \and L. Zikatanov }
\institute{Centre de Recerca Matematica (CRM), Barcelona, Spain
\texttt{bayuso@crm.cat}\\
\and Department of Mathematics, University of California at San Diego, California, USA \texttt{\{mholst, zhu\}@math.ucsd.edu}\\
\and Department of Mathematics, The Pennsylvania State University, Pennsylvania, USA \texttt{ltz@math.psu.edu}}
%
%
\maketitle

\abstract{In this paper, we present a multigrid preconditioner for solving the
  linear system arising from the piecewise linear nonconforming
  Crouzeix-Raviart discretization of second order elliptic problems
  with jump coefficients. The preconditioner uses the standard
  conforming subspaces as coarse spaces. Numerical tests show both
  robustness with respect to the jump in the coefficient and
  near-optimality with  respect to the number of degrees of freedom.
}

\section{Introduction}
\label{zhu_y_mini_3_sec:intro}

The purpose of this paper is to present a multigrid preconditioner for solving the
linear system arising from the $\mathbb{P}^{1}$ nonconforming
Crouzeix-Raviart (CR) discretization of second order elliptic problems
with jump coefficients.  The multigrid preconditioner we consider here uses pointwise
relaxation (point Gauss-Seidel/Jacobi iterative
methods) as a smoother, followed by a subspace (coarse grid) correction which uses
the standard multilevel structure for the nested $\mathbb{P}_1$
conforming finite element spaces.  The subspace correction step is
motivated by the observation that the standard $\mathbb{P}^{1}$
conforming space is a subspace of the 
CR finite element space.

One of the main benefits of this algorithm is that it is very easy to implement in practice. The procedure is the same as the standard multigrid algorithm on conforming spaces, and the only difference is the prolongation and restriction matrices on the finest level. Since the spaces are nested, the prolongation matrix is simply the matrix representation of the natural inclusion operator from the conforming space to the CR space.

The idea of using  conforming subspaces to construct preconditioners for CR discretization has been
used in \cite{Xu.J1989,Xu.J1996} in the context of smooth coefficients. For the case of jumps in the coefficients, domain decomposition preconditioners have been studied in \cite{Sarkis.M1994,sarkisNC1}  and the BPX preconditioner has been considered in \cite{Ayuso-de-Dios.B;Holst.M;Zhu.Y;Zikatanov.L2010} in connection with preconditioners for discontinuous Galerkin methods.

In the context of jump coefficients, the analysis of  multigrid preconditioners for conforming discretizations is given  in \cite{XuJ_ZhuY-2008aa}. For CR discretizations,  
the analysis is more involved due to the nonconformity of the space, and special technical tools developed in \cite{Ayuso-de-Dios.B;Holst.M;Zhu.Y;Zikatanov.L2010} are necessary. Due to space restrictions, we only state the main result (Theorem~\ref{zhu_y_mini_3_teo2} in Section~\ref{zhu_y_mini_3_sec:mg}), and provide numerical results that support it. Detailed analyses and further discussion of the algorithm will be presented in a forthcoming paper. 

The paper is organized as follows. In Section
\ref{zhu_y_mini_3_sec:cr}, we give basic notation and the finite element
discretizations. In Section \ref{zhu_y_mini_3_sec:mg}, we present the multigrid
algorithm and discuss its implementation and convergence. Finally, in
Section~\ref{zhu_y_mini_3_sec:num} we verify numerically the theoretical results by
presenting several numerical tests for two and three dimensional model
problems.

\section{Preliminaries}
\label{zhu_y_mini_3_sec:cr}
Let $\Omega\subset \mathbb{R}^{d}$ ($d=2, 3$) be an open polygonal domain. Given $f\in L^2(\Omega)$, we consider the following model problem: Find $u\in H_{0}^{1}(\Omega)$ such that
\begin{equation} \label{zhu_y_mini_3_eqn:model}
a(u, v) :=( \kappa \nabla u, \nabla v)=(f, v) \quad \forall v\in H_{0}^{1}(\Omega)\,,
\end {equation}
where the diffusion coefficient $\kappa \in L^{\infty}(\Omega)$ is assumed to be piecewise constant, namely, $\kappa(x)|_{\Omega_{m}} = \kappa_{m}$ is a constant for each (open) polygonal subdomain $\Omega_{m}$  satisfying $\cup_{m=1}^M
\overline{\Omega}_m=\overline{\Omega}$ and $\Omega_m
\cap \Omega_n =\emptyset$ for $m\neq n$. 

We assume that there is an initial (quasi-uniform) triangulation 
$\mathcal{T}_{0}$, with mesh size $h_{0}$, such that for all $T\in\mathcal{T}_{0}$ $\kappa_{T}:=\kappa(x)|_{T}$ is constant. Let $\mathcal{T}_{j}:= \mathcal{T}_{h_{j}}$ ($j =1, \cdots, J$) be a family of uniform refinement  of $\mathcal{T}_{0}$ with mesh size $h_{j}$.  Without loss of generality, we assume that the mesh size $h_{j} \simeq 2^{-j} h_{0} \;\;(j = 0, \cdots, J)$ and denote $h = h_{J}$. 

On each level $j=0, \cdots, J,$ we define $V_{j}$ as the standard $\mathbb{P}^{1}$ conforming finite element space defined on $\mathcal{T}_{j}$.  Then the standard conforming finite element discretization of \eqref{zhu_y_mini_3_eqn:model} reads: 
\begin{equation}
\label{zhu_y_mini_3_prob:c}
\mbox{Find  } u_{j}\in V_{j} \mbox{  such that  } a(u_{j}, v_{j}) = (f, v_{j}), \qquad \forall v_{j} \in V_{j}.
\end{equation}
For each $j = 0, \cdots, J$, we define the induced operator for \eqref{zhu_y_mini_3_prob:c} as 
\[
(A_j v_j,w_j) = a(v_j,w_j), \quad \forall v_j,  w_j\in V_j.
\] 

We denote $\mathcal{E}_{h}$ the set of all edges (in 2D) or faces (in 3D) of $\mathcal{T}_{h}$.
Let  $V^{CR}_{h}$ be the piecewise linear nonconforming Crouzeix-Raviart finite element space defined by: 
\begin{equation*}
V^{CR}_{h}\!=\!\left\{ v\in L^{2}(\Omega) \, : \, v_{|_{T}}
\in\, \mathbb{P}^{1}(T) \, \forall T \in \mathcal{T}_{h}\, \mbox{ and  }
\int_{e} 
\lbrack\!\lbrack{v}\rbrack\!\rbrack_{e}ds=0 \,\, \forall\, e\in \mathcal{E}_{h}\right\},
\end{equation*}
where $\mathbb{P}^{1}(T)$ denotes the space of linear polynomials on $T$ and $\lbrack\!\lbrack{v}\rbrack\!\rbrack_{e}$ denotes the jump across the edge/face $e\in \mathcal{E}_{h}$ with $\lbrack\!\lbrack{v}\rbrack\!\rbrack_{e} = v$ when $e\subset \partial \Omega$.  
In the sequel,  let us denote $V_{J+1}:= V_{h}^{CR}$ for simplicity.  We remark that all these finite element spaces are nested, that is,
$$
	V_{0} \subset\cdots\subset V_{J} \subset V_{J+1}.
$$

The $\mathbb{P}^{1}$-nonconforming finite element approximation to \eqref{zhu_y_mini_3_eqn:model} reads:
\begin{equation}\label{zhu_y_mini_3_prob:cr}
  \mbox{Find } u\in V^{CR}_{h} : a_{h}(u, w) :=\displaystyle\sum_{T\in \mathcal{T}_{J}} \int _{T}\kappa_{T}\nabla u \cdot \nabla w=(f, w),  \forall\, w\in V^{CR}_{h}. 
 \end{equation}
The bilinear form $a_{h}(\cdot, \cdot)$ induced a natural energy norm: $|v|_{h,\kappa}:=\sqrt{a_{h}(v,v)}$ for any $v\in V^{CR}_{h}$. 
In operator form, we are going to solve the linear system 
\begin{equation}
\label{zhu_y_mini_3_eqn:op}
	Au = f,
\end{equation}
where $A$ is the operator induced by \eqref{zhu_y_mini_3_prob:cr}, namely
$$
	(A v, w) = a_{h}(v ,w),\qquad \forall v, w\in V_{h}^{CR}.
$$
\section{A Multigrid Preconditioner}
\label{zhu_y_mini_3_sec:mg}
The action of the standard multigrid $V$-cycle preconditioner
$B:=B_{J+1}:V_{J+1} \mapsto V_{J+1}$ on a given $g\in V_{J+1}$ is
recursively defined by the following algorithm (cf. \cite{Bramble.J1993}):
\begin{algorithm}[$V$-cycle]
\label{zhu_y_mini_3_alg::vcycle} 
Let $g_{J+1}=g$, and $B_0=A_0^{-1}.$ For $j = 1, \cdots, J+1,$ we define
recursively $B_j g_j$ for any $g_j\in V_{j}$ by the following three steps:
\begin{enumerate}
  \item Pre-smoothing : $w_1=R_j g_j;$
  \item Subspace correction: $w_2=w_1+ B_{j-1}Q_{j-1}(g_j - A_j w_1);$
  \item Post-smoothing: $B_jg_j:=w_2+R_j^*(g_j-A_j w_2).$
\end{enumerate}
\end{algorithm}
In this algorithm, $R_{j}$ corresponds to a Gauss-Seidel or
a Jacobi iterative method known as a smoother; and $Q_{j}$ is the
standard $L^{2}$ projection on $V_j$:
$$
	(Q_{j} v, w_{j} ) = (v, w_{j}), \qquad \forall w_{j} \in V_{j}, \;\; (j=0, \cdots, J).
$$

The implementation of Algorithm~\ref{zhu_y_mini_3_alg::vcycle} is almost identical to the implementation of the standard multigrid $V$-cycle (cf. \cite{Briggs.W;Henson.V;McCormick.S2000}). Between the conforming spaces, we use the standard prolongation and restriction matrices (for conforming finite elements). The corresponding matrices between $V_{J}$ and $V_{J+1}$, are however different.
The prolongation matrix on $V_{J}$ can be viewed as the matrix representation of the natural inclusion $\mathcal{I}_{J}: V_{J} \to V_{J+1},$ which is defined by
$$(\mathcal{I}_{J} v)(x) = \sum_{e\in \mathcal{E}_{h}} v(m_{e}) \psi_{e}(x),$$
where $\psi_{e}$ is the CR basis on the edge/face $e\in \mathcal{E}_{h}$ and $m_{e}$ is the barycenter of $e$. Therefore, the prolongation matrix has the same sparsity pattern as the edge-to-vertex (in 2D), or face-to-vertex (in 3D) connectivity, and each nonzero entry in this matrix equals the constant $1/d$ where $d$ is the space dimension. The restriction matrix is simply the transpose  of the prolongation matrix. 

The efficiency and robustness of this preconditioner can be  analyzed in terms of the \emph{effective condition number} (cf. \cite{XuJ_ZhuY-2008aa}) defined as follows:
\begin{definition}
	\label{zhu_y_mini_3_def::eff}
	Let $V$ be a real $N$ dimensional Hilbert space, and $S:V\to V$ be a symmetric positive definition operator with eigenvalues $0 < \lambda_{1} \le \cdots \le \lambda_{N}.$ The $m$-th \emph{effective condition} number of $S$ is defined by 
	$$ \mathcal{K}_{m}(S):=\lambda_{N}(S)/\lambda_{m+1}(S).$$
\end{definition}
Note that the standard condition number $\mathcal{K}(BA)$ of the preconditioned system $BA$ will be large due to the large jump in the coefficient $\kappa$. However, there might be only a small (fixed) number of small eigenvalues of $BA$, which cause the large condition number; and the other eigenvalues are bounded nearly uniformly. In particular, we have the following main result:
\begin{theorem}\label{zhu_y_mini_3_teo2}
  Let $B$ be the multigrid $V$-cycle preconditioner defined in
  Algorithm~\ref{zhu_y_mini_3_alg::vcycle}. Then there exists a fixed integer
  $m_0<M,$ depending only on the distribution of the coefficient
  $\kappa$, such that 
	\begin{equation*}
		\mathcal{K}_{m_{0}}(B A) \le C^{2} |\log h|^{2} = C^{2} J^{2}\;,
	\end{equation*}
	where the constant $C>0$ is independent of the coefficients and mesh size.
\end{theorem}
The analysis is based on the subspace correction framework \cite{Xu.J1992a},  but some technical tools developed in \cite{Ayuso-de-Dios.B;Holst.M;Zhu.Y;Zikatanov.L2010} are needed to deal with nonconformity of the finite element spaces. Due to space restriction, a detailed analysis will be reported somewhere else.

Thanks to Theorem~\ref{zhu_y_mini_3_teo2} and a standard PCG convergence result (cf. \cite[Section 13.2]{Axelsson.O1994}), the PCG algorithm with the multigrid $V$-cycle preconditioner defined in Algorithm~\ref{zhu_y_mini_3_alg::vcycle} has the following convergence estimate:
$$|u-u_i|_{h,\kappa}\le 2(\mathcal{K}(BA)-1)^{m_0}
\left(\frac{CJ-1}{CJ+1}\right)^{i-{m_0}}|u-u_0|_{h,\kappa}\;,
$$
where $u_{0}$ is the initial guess, and $u_{i}$ is the solution of
$i$-th PCG iteration.  Although the condition number $\mathcal{K}(BA)$ might be large, the convergence rate of the PCG
algorithm is asymptotically dominated by $\frac{CJ-1}{CJ+1},$ which is
determined by the effective condition number
$\mathcal{K}_{m_{0}}(BA)$. Moreover, this bound of asymptotic convergence rate convergence is independent of the coefficient $\kappa$, but depends on the mesh size logarithmically.

\section{Numerical Results}
\label{zhu_y_mini_3_sec:num}
In this section, we present several numerical tests in 2D and 3D which verify the result in Theorem~\ref{zhu_y_mini_3_teo2} on the performance of the
multigrid $V$-cycle preconditioner described in the previous
sections. The numerical tests show that the effective condition
numbers of the preconditioned linear systems (with $V$-cycle
preconditioner) are nearly uniformly bounded.
\subsection{A 2D Example}
\begin{figure}[htbp]
\centering
	\parbox{0.45\textwidth}{
       \includegraphics[width=0.47\textwidth]{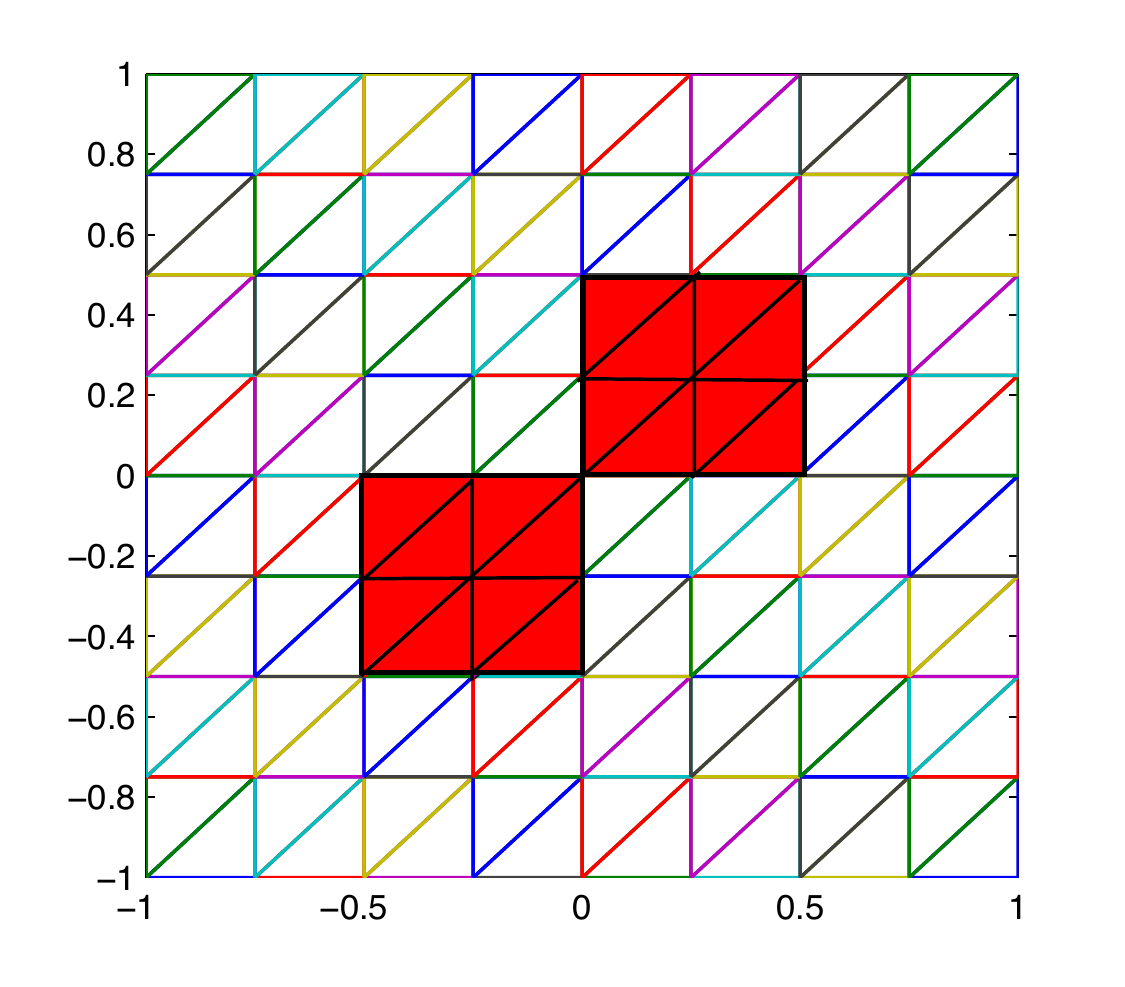}
       \caption{2D Computational Domain}
       \label{zhu_y_mini_3_fig:domain2d}}
       \quad
       \begin{minipage}{0.45\textwidth}
       \includegraphics[width=0.99\textwidth]{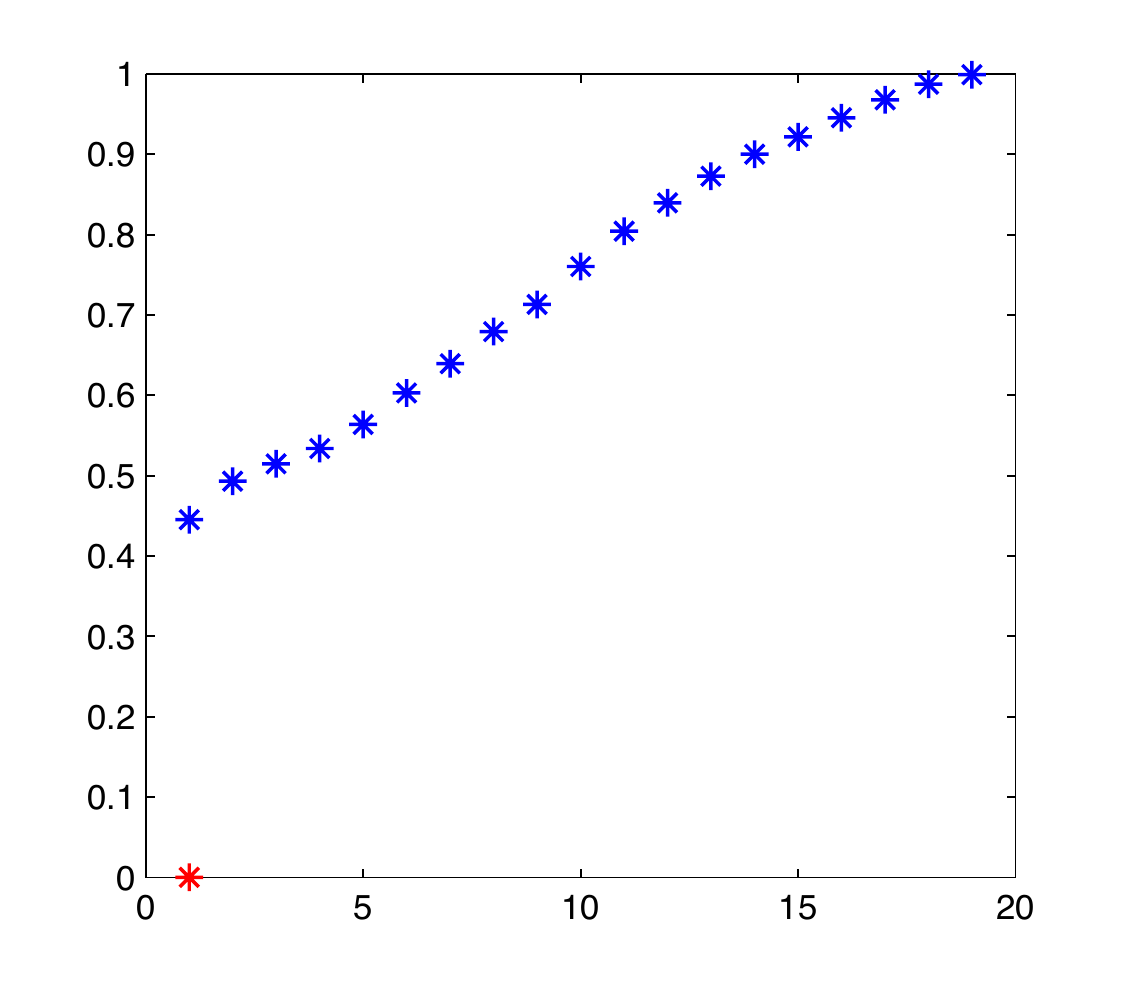}
       \caption{Eigenvalue Distribution of $BA$}
       \label{zhu_y_mini_3_fig:eig2d}
       \end{minipage}
\end{figure}

As a first model problem, we consider equation~\eqref{zhu_y_mini_3_eqn:model} in the square
$\Omega = (-1,1)^{2}$ with coefficient such that, $\kappa(x) =1$ for
$x\in \Omega_1=(-0.5,0)^{2}\cup (0, 0.5)^{2}$, and $\kappa(x) = \epsilon$ for
$x$ in the remaining subdomain, $x\in \Omega\setminus\Omega_1$ (see Figure~\ref{zhu_y_mini_3_fig:domain2d}). By
decreasing the value of $\epsilon$ we increase the contrast in the PDE
coefficients.  

Our initial triangulation on level 0 has mesh size $h_{0} =
2^{-1}$ and resolves the interfaces where the coefficients have
discontinuities. Then on
each level, we uniformly refine the mesh by subdividing each element
into four congruent children. In this example, we use 1 forward/backward
Gauss-Seidel iteration as pre/post smoother in the multigrid
preconditioner, and the stopping criteria of the PCG algorithm is
$\|r_{k}\| / \|r_{0}\| <10^{-7}$ where $r_{k}$ is the the residual at $k$-th
iteration.

Figure~\ref{zhu_y_mini_3_fig:eig2d} shows the eigenvalue distribution of the
multigrid $V$-cycle preconditioned system $BA$ when $h=2^{-5}$ (level =4) and
$\epsilon = 10^{-5}.$ As we can see from this figure, there is only
one small eigenvalue that deteriorates with respect to the jump in the coefficient and the mesh size.

Table \ref{zhu_y_mini_3_tab:1} shows the estimated condition number $\mathcal{K}$
and the effective condition number $\mathcal{K}_{1}$ of $BA$. It can be observed that the condition number $\mathcal{K}$ increases rapidly with respect to the increase of the jump in the coefficients and the number of degrees of freedom. On the other hand, the number of PCG iterations increases only a small amount,  and the corresponding effective condition number is
nearly uniformly bounded, as predicted by Theorem~\ref{zhu_y_mini_3_teo2}. 
\begin{table}
	\begin{tabular}{c|c||c|c|c|c|c}
\hline
{$\epsilon$} & levels &  0         & 1        & 2        & 3        & 4 \\
\hline\hline
\multirow{2}{*}{$1$}
 & $\mathcal{K}$ &  1.65 (8)&  1.83 (10)&   1.9 (10)&   1.9 (10)&  1.89 (10)\\
 & $\mathcal{K}_{1}$ &  1.44 &  1.78 &  1.77 &  1.78 &  1.76 \\\hline
\multirow{2}{*}{$10^{-1}$}
 & $\mathcal{K}$ &  3.78 (10)&  3.69 (11)&  3.76 (12)&  3.79 (12)&  3.88 (12)\\
 & $\mathcal{K}_{1}$ &  1.89 &  1.87 &  1.93 &  1.92 &  1.95 \\\hline
\multirow{2}{*}{$10^{-2}$}
 & $\mathcal{K}$ &  23.4 (12)&  23.6 (13)&  24.6 (13)&  25.1 (14)&    26 (15)\\
 & $\mathcal{K}_{1}$ &  2.15 &  1.96 &  1.99 &  1.97 &  2.24 \\\hline
\multirow{2}{*}{$10^{-3}$}
 & $\mathcal{K}$ &   218 (13)&   223 (14)&   232 (15)&   238 (16)&   246 (16)\\
 & $\mathcal{K}_{1}$ &  2.19 &  1.98 &     2 &  1.98 &  2.29 \\\hline
\multirow{2}{*}{$10^{-4}$}
 & $\mathcal{K}$ & 2.17e+03 (14)& 2.21e+03 (15)& 2.31e+03 (16)& 2.37e+03 (18)& 2.45e+03 (18)\\
 & $\mathcal{K}_{1}$ &   2.2 &  1.98 &     2 &  1.98 &   2.3 \\\hline 
\multirow{2}{*}{$10^{-5}$}
 & $\mathcal{K}$ & 2.17e+04 (15)& 2.21e+04 (16)& 2.31e+04 (17)& 2.37e+04 (19)& 2.76e+04 (19)\\
 & $\mathcal{K}_{1}$ &   2.2 &  1.98 &     2 &  1.98 &  2.64 \\\hline
\end{tabular}\label{zhu_y_mini_3_tab:1}
\caption{Estimated condition number $\mathcal{K}$ (number of PCG iterations) and the effective condition number $\mathcal{K}_{1}$}
\end{table}

\subsection{A 3D Example}
\begin{figure}[htbp]
\centering
	\parbox{0.45\textwidth}{
       \includegraphics[width=0.43\textwidth]{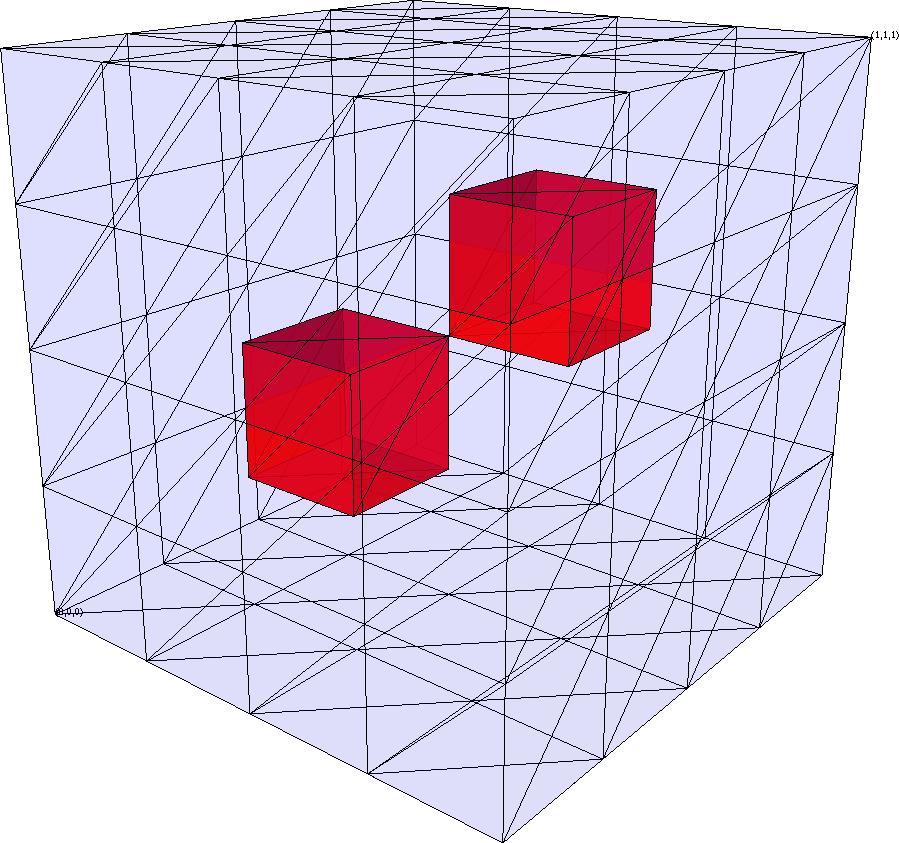}
       \caption{3D Computational Domain}
       \label{zhu_y_mini_3_fig:domain3d}}
       \quad
       \begin{minipage}{0.45\textwidth}
       \includegraphics[width=0.99\textwidth]{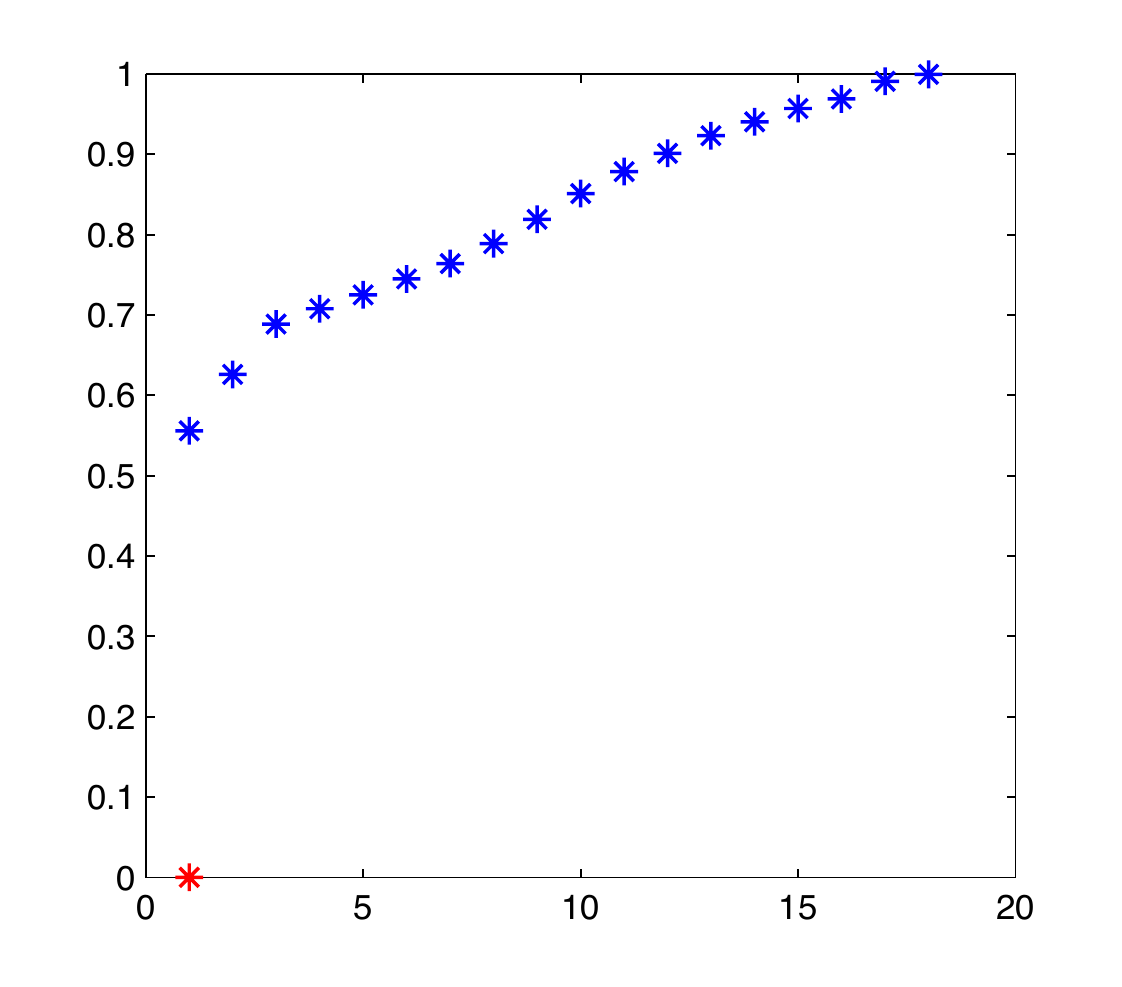}
       \caption{Eigenvalue Distribution of $BA$}
       \label{zhu_y_mini_3_fig:eig3d}
       \end{minipage}
\end{figure}

In this second example, we consider the model problem~\eqref{zhu_y_mini_3_eqn:model} in
the open unit cube in 3D with a similar setting for the coefficient.  We set
$\kappa(x) =1$ for $x\in \Omega_{1} =(0.25, 0.5)^{3}$ or $x\in
\Omega_{2}=(0.5, 0.75)^{3}$, and $\kappa(x) =\epsilon$ for the
remaining subdomain (that is, for $x\in
\Omega\setminus(\Omega_{1}\cup\Omega_2)$). The domain $\Omega$ and the subdomains just
described are shown in Figure~\ref{zhu_y_mini_3_fig:domain3d}. The coarsest partition
has mesh size $h_{0} =2^{-2}$, and it is set in a way so that it resolves the
interfaces where the coefficient has jumps. 

To test the effects of the smoother, in this example we used 5
forward/backward Gauss-Seidel as smoother in the multigrid
preconditioner. In order to test more severe jumps in the coefficients, we set the stopping criteria $\|r_{k}\| / \|r_{0}\|<10^{-12}$ for the PCG algorithm in this experiment.

Figure~\ref{zhu_y_mini_3_fig:eig3d} shows the eigenvalue distribution of the multigrid $V$-cycle preconditioned system $BA$ when $h=2^{-5}$ (level=3) and $\epsilon = 10^{-5}.$ As before, this figure shows that there is only one small eigenvalue that even deteriorates with respect to the jump in the coefficients and the mesh size. 

\begin{table}
\centering
\begin{tabular}{c|c||c|c|c|c}
\hline
 $\epsilon$ & levels &  0         & 1        & 2        & 3        \\
\hline\hline
\multirow{2}{*}{$1$}
 & $\mathcal{K}$ & 1.19 (8)& 1.34 (11)& 1.37 (11)& 1.36 (11)\\
 &$\mathcal{K}_{1}$ & 1.16 &  1.26 &  1.31& 1.29\\\hline
\multirow{2}{*}{$10^{-1}$}
 & $\mathcal{K}$ & 2.3 (10)& 1.94(13)& 1.75 (13)& 1.67 (14)\\
 & $\mathcal{K}_{1}$ & 1.60 & 1.56 & 1.45 &  1.43 \\\hline
\multirow{2}{*}{$10^{-3}$}
 & $\mathcal{K}$ &   86.01 (11)&   63.07 (16)&  52.67 (17)& 48.19(17)\\
 & $\mathcal{K}_{1}$ & 2.4 &  2.12  &   1.89 &  1.78 \\\hline
\multirow{2}{*}{$10^{-5}$}
 & $\mathcal{K}$ &  8.39+03 (13)&  6.15e+03 (18)&  5.13e+03 (19)&   4.70e+03(19)\\
 &  $\mathcal{K}_{1}$ &   2.44 &   2.14 &  1.91 & 1.80 \\\hline
\multirow{2}{*}{$10^{-7}$}
 & $\mathcal{K}$ &  8.39+05 (14)&  6.15e+05 (21)&  5.13e+05 (23)&  4.70e+05(21)\\
 & $\mathcal{K}_{1}$ &  2.45 &  2.14  & 1.91  & 1.80\\\hline
\end{tabular}
\caption{Estimated condition number $\mathcal{K}$ (number of PCG iterations) and effective condition number $\mathcal{K}_{1}$.}
\label{zhu_y_mini_3_tab:3d}
\end{table}
Table \ref{zhu_y_mini_3_tab:3d} shows the estimated condition number $\mathcal{K}$
(with the number of PCG iterations), and the effective condition
number $\mathcal{K}_{1}$. As is easily seen from the results in
this table, the condition number $\mathcal{K}$ increases when
$\epsilon$ decreases, i.e. the condition number grows when the
jump in the coefficients becomes larger. On the other hand, the
results in Table~\ref{zhu_y_mini_3_tab:3d} show that the effective condition number
$\mathcal{K}_{1}$ remains nearly uniformly bounded with respect to the
mesh size and it is robust with respect to the jump in the coefficient, thus
confirming the result stated in Theorem~\ref{zhu_y_mini_3_teo2}: a PCG with multigrid
$V$-cycle preconditioner provides a robust, nearly optimal solver for the CR 
approximation to \eqref{zhu_y_mini_3_prob:cr}.

\subsection*{Acknowledgments}
\vskip -0.2cm
First author has been supported by MEC grant MTM2008-03541 and 2009-SGR-345 
from AGAUR-Generalitat de Catalunya. 
The work of the second and third authors was supported in part by
NSF/DMS Awards 0715146 and 0915220, and by DOD/DTRA Award HDTRA-09-1-0036.
The work of the fourth author was supported in part by the
NSF/DMS Award 0810982.

%
\bibliographystyle{plainnat}
\bibliography{zhu_y_mini_3} 

\begin{thebibliography}{10}
\providecommand{\natexlab}[1]{#1}
\providecommand{\url}[1]{\texttt{#1}}
\expandafter\ifx\csname urlstyle\endcsname\relax
  \providecommand{\doi}[1]{doi: #1}\else
  \providecommand{\doi}{doi: \begingroup \urlstyle{rm}\Url}\fi

\bibitem[Axelsson(1994)]{Axelsson.O1994}
O.~Axelsson.
\newblock \emph{Iterative solution methods}.
\newblock Cambridge University Press, Cambridge, 1994.
\newblock ISBN 0-521-44524-8.

\bibitem[Ayuso~de Dios et~al.(2010)Ayuso~de Dios, Holst, Zhu, and
  Zikatanov]{Ayuso-de-Dios.B;Holst.M;Zhu.Y;Zikatanov.L2010}
B.~Ayuso~de Dios, M.~Holst, Y.~Zhu, and L.~Zikatanov.
\newblock {Multilevel Preconditioners for Discontinuous Galerkin Approximations
  of Elliptic Problems with Jump Coefficients}.
\newblock \emph{Arxiv preprint arXiv:1012.1287}, 2010.

\bibitem[Bramble(1993)]{Bramble.J1993}
J.~H. Bramble.
\newblock \emph{Multigrid Methods}, volume 294 of \emph{Pitman Research Notes
  in Mathematical Sciences}.
\newblock Longman Scientific \& Technical, Essex, England, 1993.

\bibitem[Briggs et~al.(2000)Briggs, Henson, and
  McCormick]{Briggs.W;Henson.V;McCormick.S2000}
W.~L. Briggs, V.~E. Henson, and S.~F. McCormick.
\newblock \emph{A multigrid tutorial}.
\newblock Society for Industrial and Applied Mathematics (SIAM), Philadelphia,
  PA, second edition, 2000.
\newblock ISBN 0-89871-462-1.

\bibitem[Sarkis(1994{\natexlab{a}})]{sarkisNC1}
M.~Sarkis.
\newblock Multilevel methods for {$P_1$} nonconforming finite elements and
  discontinuous coefficients in three dimensions.
\newblock In \emph{Domain decomposition methods in scientific and engineering
  computing ({U}niversity {P}ark, {PA}, 1993)}, volume 180 of \emph{Contemp.
  Math.}, pages 119--124. Amer. Math. Soc., Providence, RI, 1994{\natexlab{a}}.

\bibitem[Sarkis(1994{\natexlab{b}})]{Sarkis.M1994}
M.~V. Sarkis.
\newblock \emph{Schwarz Preconditioners for Elliptic Problems with
  Discontinuous Coefficients Using Conforming and Non-Conforming Elements}.
\newblock PhD thesis, Courant Institute of Mathematical Science of New York
  University, 1994{\natexlab{b}}.

\bibitem[Xu(1989)]{Xu.J1989}
J.~Xu.
\newblock \emph{Theory of Multilevel Methods}.
\newblock PhD thesis, Cornell University, 1989.

\bibitem[Xu(1992)]{Xu.J1992a}
J.~Xu.
\newblock Iterative methods by space decomposition and subspace correction.
\newblock \emph{SIAM Review}, 34:\penalty0 581--613, 1992.

\bibitem[Xu(1996)]{Xu.J1996}
J.~Xu.
\newblock The auxiliary space method and optimal multigrid preconditioning
  techniques for unstructured meshes.
\newblock \emph{Computing}, 56:\penalty0 215--235, 1996.

\bibitem[Xu and Zhu(2008)]{XuJ_ZhuY-2008aa}
J.~Xu and Y.~Zhu.
\newblock Uniform convergent multigrid methods for elliptic problems with
  strongly discontinuous coefficients.
\newblock \emph{Math. Models Methods Appl. Sci.}, 18\penalty0 (1):\penalty0 77
  --105, 2008.

\end{thebibliography}

\end{document}